\documentclass[11pt,draft]{amsart}

\usepackage[cp1251]{inputenc}
\usepackage[ukrainian,english]{babel}
\usepackage{amsmath,amsthm,amsfonts,amssymb}

\theoremstyle{plain}
\newtheorem{theorem}{Theorem}
\newtheorem{lemma}{Lemma}
\newtheorem{corollary}{Corollary}
\theoremstyle{definition}
\newtheorem{definition}{Definition}

\newtheorem{property}{}
\theoremstyle{remark}

\newcommand{\N}{\mathbb{N}}
\newcommand{\Z}{\mathbb{Z}}
\newcommand{\set}[1]{\left\{#1\right\}}
\newcommand{\abs}[1]{\left\lvert#1\right\rvert}
\newcommand{\cl}{\mathrm{cl}}
\renewcommand{\O}[1]{\mathrm{O}^{#1}}
\newcommand{\bO}[1]{\bar{\mathrm{O}}^{#1}}
\newcommand{\bcyl}[2]{\bar{\mathrm{O}}^{#1}_{[#2]}}
\newcommand{\bint}[2]{\bar{\mathrm{O}}^{#1}_{(#2)}}
\makeatletter
\newcommand{\Cset}[2][\bO1]{C
  [#1,
  \if\relax\expandafter\@gobble#2\relax #2\else\{#2\}\fi
  ]}%
\makeatother

\begin{document}

\title[The Ostrogradsky series and related probability measures]{The Ostrogradsky series\\ and related probability measures}
\author[S.Albeverio, O.Baranovskyi, M.Pratsiovytyi, G.Torbin]
{Sergio Albeverio$^{1,2,3,4}$, Oleksandr Baranovskyi$^{5,6}$,\\
Mykola~Pratsiovytyi$^{7,8}$ \and Grygoriy Torbin$^{1,9,10}$}

\begin{abstract}
We develop a metric and probabilistic theory for the  Ostrogradsky
representation of real numbers, i.e., the expansion of a real
number $x$ in the following form:
\begin{align*}
x&= \sum_n\frac{(-1)^{n-1}}{q_1q_2\dots q_n}=\\
&=\sum_n\frac{(-1)^{n-1}}{g_1(g_1+g_2)\dots(g_1+g_2+\dots+g_n)}\equiv
\bO1(g_1,g_2,\dots,g_n,\dots),
\end{align*}
where $q_{n+1}>q_n\in\N$, $g_1=q_1$, $g_{k+1}=q_{k+1}-q_k$. We
compare this representation with the corresponding one in terms of
continued fractions.

 We establish basic metric relations (equalities
and inequalities for ratios  of the length of cylindrical sets). We
also compute the Lebesgue measure of subsets belonging to  some
classes of closed nowhere dense sets defined by characteristic
properties of the $\bO1$-representation. In particular, the
conditions for the set $\Cset{V}$, consisting of real numbers whose
$\bO1$-symbols take values from the set $V \subset N$, to be of zero
resp. positive Lebesgue measure are found. For a random variable
$\xi$ with independent  $\bO1$-symbols $g_n(\xi)$ we prove the
theorem establishing the purity of the distribution. In the case of
singularity the conditions for such distributions to be of Cantor
type are also found.
\end{abstract}

\maketitle

$^1$~Institut f\"{u}r Angewandte Mathematik, Universit\"{a}t Bonn,
Wegelerstr. 6, D-53115 Bonn (Germany); $^2$~SFB 611, Bonn, BiBoS,
Bielefeld - Bonn; $^3$~CERFIM, Locarno and Acc. Arch., USI
(Switzerland); $^4$~IZKS, Bonn; E-mail: albeverio@uni-bonn.de

$^{5}$ ~Institute for Mathematics of NASU, Tereshchenkivs'ka str. 3,
01601 Kyiv (Ukraine); $^6$~National Pedagogical University, Pyrogova
str. 9, 01030 Kyiv (Ukraine); E-mail: ombaranovskyi@ukr.net

$^{7}$~National Pedagogical University, Pyrogova str. 9, 01030 Kyiv
(Ukraine); $^8$~~Institute for Mathematics of NASU,
Tereshchenkivs'ka str. 3, 01601 Kyiv (Ukraine); E-mail: m\_
pratz@ukr.net

$^9$~National Pedagogical University, Pyrogova str. 9, 01030 Kyiv
(Ukraine) $^{10}$~Institute for Mathematics of NASU,
Tereshchenkivs'ka str. 3, 01601 Kyiv (Ukraine); E-mail:
torbin@wiener.iam.uni-bonn.de (corresponding author)

\medskip


\medskip


\section*{Introduction}

There are many different methods for the expansions  and encodings
(representations) of real numbers by using a finite as well as
infinite alphabet $A$. The $s$-adic expansions,
continued fractions, $f$\nobreakdash-expansions, the L\"{u}roth
expansions  etc. are widely used in mathematics (see, e.g.,
\cite{Sch95}). Each representation has its own specificity and it
generates its ``own geometry'' and metric theory. To each
representation there is associated system of cylindrical sets,
which forms a system of partitions of the unit interval (real
line). We also have a corresponding ``coordinate system'' (a union
of conditions for the determination of the position of a point)
which is a convenient tool for the description of a wide class of
fractals in a simple formal way. From the ratios of the lengths of
cylindrical sets  the basic metric relations follow (in the form
of equalities and inequalities) which are crucial for the
development of the corresponding metric theory, i.e. a theory
about measure (e.g., Jordan, Lebesgue, Hausdorff,
Hausdorff-Billingsley,...) of sets of real numbers defined by
characteristic properties of their digits in  the corresponding
representation.

Let $A$ be an alphabet of symbols for the representation of
numbers in some fixed system of representation, let $\alpha_k(x)$
be the $k$-th symbol of the representation
$\Delta_{\alpha_1...\alpha_k...}$ of a real number $x\in[0,1]$,
let $N_i(x,n)$ be the number of the symbol ``$i$'' among the first
$n $ symbols of the representation of $x$, $i\in A$, and let
$\nu_i(x)=\lim\limits_{n\rightarrow\infty} \frac {N_i(x,n)}{n}$
(supposed to exist).

The following sets are traditional objects for the investigations
in the metric theory:
\begin{align*}
&E_f\left(\genfrac{}{}{0pt}{}{k_1k_2\dots k_n}{c_1c_2\dots c_n}
\right)= \set{x : \alpha_{k_i}(x)=c_i\in A},\\
&C[f,\{V_k\}]=\set{x : \alpha_k(x)\in V_k\subset A},\\
&M[f,\tau]=\set{x : \nu_i(x)=\tau_i,\; \forall i\in A;\; \tau =
(\tau_1, \tau_2, \dots ), \tau_i \geq 0, \sum\limits_i \tau_i = 1},\\
 &T[f,\bar\nu]=\set{x : \text{$\nu_i(x)$
does not exist for all symbols from $A$}}.
\end{align*}

Let us remark that $f$ denotes the method of the representation of
numbers ($f: L\rightarrow[0,1]$, where $L=A\times A\times\dotsb$).
During recent years the  interest into the latter set has been
considerably  increasing (see, e.g., \cite{APT3,APTUMJ, Ols04b,
Ols04a,PT95}).

The presented paper devoted to the investigation of the expansion
of real numbers in the first Ostrogradsky series (they were
introduced by M.~V.~Ostrogradsky, a well known ukrainian
mathematician  who lived from 1801 to 1862). We shall also present
the development of the corresponding metric and probabilistic
theory. In this case the alphabet $A$ coincides with the set $\N$
of positive integers.

The expansion of $x$ of the form:
\begin{equation}\label{?}
x=\frac{1}{q_1}-\frac{1}{q_1q_2}+\dots
+\frac{(-1)^{n-1}}{q_1q_2\dots q_n}+\dotsb,
\end{equation}
where $q_n$ are positive integers and $q_{n+1}>q_n$ for all $n$,
is said to be the expansion of $x$ in the first Ostrogradsky
series. The expansion of $x$ of the form:
\begin{equation}\label{??}
x=\frac{1}{q_1}-\frac{1}{q_2}+\dots
+\frac{(-1)^{n-1}}{q_n}+\dotsb,
\end{equation}
where $q_n$ are positive integers and $q_{n+1}\geq q_n(q_n+1)$ for
all $n$, is said to be the expansion of $x$ in the second
Ostrogradsky series. Each irrational number has a unique expansion
of the form~\eqref{?} or~\eqref{??}. Rational numbers have two
finite different representations of the above form (see, e.g.,
\cite{Rem51a}).

Equality \eqref{?} can be rewritten in the following way:
\begin{equation}\label{???}
x=\frac1{g_1}-\frac1{g_1(g_1+g_2)}+\dots+\frac{(-1)^{n-1}}{g_1(g_1+g_2)\dots(g_1+g_2+\dots+g_n)}+\dotsb,
\end{equation}
where $g_1=q_1$, $g_{n+1} = q_{n+1}-q_n$ for any $n\in\N$. The
expression~\eqref{???} is said to be the $\bO1$-representation and
the symbol $g_n=g_n(x)$ is said to be the $n$-th $\bO1$-symbol of
$x$.

Shortly before his death, M.~Ostrogradsky has proposed two
algorithms for the representation of real numbers via  alternating
series of the form \eqref{?} and \eqref{??}, but he did not publish
any papers on this problems. Short Ostrogradsky's remarks concerning
the above representations have been found by E.~Ya.~Remez
\cite{Rem51a} in the hand-written fund of the Academy of Sciences of
USSR. E.~Ya.~Remez has pointed out some similarities between the
Ostrogradsky series and continued fractions. He also paid a great
attention to the applications of the Ostrogradsky series for the
numerical methods for solving algebraic equations. In the editorial
comments to the book \cite{Khi97} B.~Gnedenko has pointed out that
there are no fundamental investigations of properties of the above
mentioned representations. Analogous problems were studied by
W.~Sierpi\'{n}ski \cite{Sie74} and T.~A.~Pierce \cite{Pie29}
independently. Some algorithms for the representation of real
numbers in positive and alternating series were proposed in
\cite{Sie74}. Two of these algorithms lead to the Ostrogradsky
series \eqref{?} and \eqref{??}. An algorithm also leading to the
representation of irrational numbers in the form of the series
\eqref{?} has been considered in \cite{Pie29}.

There exists a series of papers devoted to the applications of the
Ostrogradsky series. Let us mention some of them. Connections
between the Ostrogradsky algorithms and the algorithm for the
continued fractions have been established in \cite{BS74}. This book
contains also generalizations of the above algorithms. In the paper
\cite{Mel74} different types of $p$-adic continued fractions have
been constructed on the basis of $p$-adic analogs of Euclid and
Ostrogradsky algorithms. Combining in a special way the algorithms
of Engel and Ostrogradsky, the same author in the paper \cite{Mel76}
has constructed an algorithm for the representation of real numbers
via series which converge faster then the corresponding Engel's and
Ostrogradsky's series.  \cite{VZ75} is devoted to the investigation
of the first Ostrogradsky algorithm and to the determination of the
expectation of the random variables $(q_j+1)^\nu$, $\nu\geq 0$ and
$r_n=\sum\limits_{j=n+1}^\infty \frac{(-1)^{j+1}}{q_1q_2\dots q_j}$,
where $q_j = q_j(\alpha)$ are  random variables depending on the
random variable $\alpha$, uniformly distributed on the unit
interval.
 In the same paper a
generalization of the Ostrogradsky algorithms for approximations in
Banach spaces has been proposed.

 In the presented paper we study basic metric relations (equalities and inequalities
for ratios  of the length of corresponding cylinders) for the
$\bO1$-representation of reals. In Section 2 we paid the main
attention to the problem of the approximation of real numbers by
partial sums  of the Ostrogradsky series. We stress some
similarities of the $\bO1$-representation with the continued
fraction representation. Recurrent formulas for the
$\bO1$-convergents (analogs of the convergents for continued
fractions) is also studied in this Section.

In Section 3 we prove basic metric relations of the
$\bO1$-representation and compare them with the corresponding
relations for continued fractions.

Sections 4 and 5 are the main ones of the paper. Section 4 is
devoted to the study of the set $\Cset{V_k}$, consisting of real
numbers whose k-th $\bO1$-symbols take values from the set
$V_k\subset N$. The central object of this section is the set
$C=\Cset{V}$ which is a particular case of the previous one (for
$V_n=V, \forall n\in\N$). Conditions for the set $\Cset{V}$ to be
of zero resp. positive Lebesgue measure $\lambda$ are found. In
particular, we prove that $\lambda(C)>0$, if $V=\{m+1, m+2,
\dots\}$,
 where $m$ is an arbitrary
positive integer. This fact stresses an essential difference between
the metric theories of continued fractions and
$\bO1$-representations.

In Section 5 we study the random variable $\xi$ with independent
$\bO1$-symbols $\xi_{k}$. In particular, we  prove that the random
variable $\xi$ with independent $\bO1$-symbols is of pure type,
i.e., it is either pure singular continuous, or pure absolutely
continuous or pure atomic. On the basis of results of the previous
Sections we study properties of the topological support of the
random variable $\xi$. In the atomic case we completely describe the
set of all atoms of the distribution. In the continuous case we give
sufficient conditions for $\xi$ to be a singular continuous
distribution of the Cantor type.

\section{Representations of real numbers by the Ostrogradsky series}

\begin{definition}
A finite or an infinite expression
\begin{equation}\label{eq:o1series}
\sum_n\frac{(-1)^{n-1}}{q_1q_2\dots q_n} =
\frac{1}{q_1}-\frac{1}{q_1q_2}+\dotsb,
\end{equation}
where $q_n$ are natural and $q_{n+1}>q_n$ for all $n$, is called
\emph{the first Ostrogradsky series} (in the sequel \emph{the
Ostrogradsky series}). The numbers $q_n$ are called \emph{the
symbols of the Ostrogradsky series} \eqref{eq:o1series}.
\end{definition}

We denote the expression \eqref{eq:o1series} briefly by
\[\O1(q_1,q_2,\dots,q_n)\]
if it contains a finite number of terms, and we speak in this case
of a finite Ostrogradsky series. We denote \eqref{eq:o1series} by
\[\O1(q_1,q_2,\dots,q_n,\dots)\]
if it contains an infinite number of terms.

Every Ostrogradsky series is convergent and its sum belongs to
$[0,1]$.

\begin{theorem}[\cite{Rem51a}]
Any real number $x\in(0,1)$ can be represented in the form
\eqref{eq:o1series}. If $x$ is irrational then the expression
\eqref{eq:o1series} is unique and it has an infinite number of
terms. If $x$ is rational then it can be represented in the form
\eqref{eq:o1series} in the following different ways:
\[x=\O1(q_1,q_2,\dots,q_{n-1}, q_{n}, q_{n}+1)=\O1(q_1,q_2,\dots,q_{n-1},q_n+1).\]
\end{theorem}

We can find the symbols of the Ostrogradsky series for a given
number $x$ using the following algorithm:
\begin{equation}\label{eq:o1alg}
\begin{array}{lll}
1=q_1x+\alpha_1            && \left(0\leq\alpha_1<x\right),\\
1=q_2\alpha_1+\alpha_2     && \left(0\leq\alpha_2<\alpha_1\right),\\
\hdotsfor{3}\\
1=q_n\alpha_{n-1}+\alpha_n && \left(0\leq\alpha_n<\alpha_{n-1}\right),\\
\hdotsfor{3}
\end{array}
\end{equation}

Let
\[
g_1=q_1\quad \text{and}\quad g_{n+1}=q_{n+1}-q_n \quad \text{for
any $n\in\N$}.
\]
Then one can rewrite series \eqref{eq:o1series} in the form
\begin{equation}\label{eq:bo1series}
\sum_n\frac{(-1)^{n-1}}{g_1(g_1+g_2)\dots(g_1+g_2+\dots+g_n)}=
\frac1{g_1}-\frac1{g_1(g_1+g_2)}+ \dotsb.
\end{equation}
We denote the expression \eqref{eq:bo1series} by
\[\bO1(g_1,g_2,\dots,g_n,\dots).\]
A representation of a number $x\in(0,1)$ by expression
\eqref{eq:bo1series} is called \emph{the $\bO1$-representation}.
The number $g_n=g_n(x)$ is called \emph{$n$-th $\bO1$-symbol} of
the number $x$.

\textbf{Definition.} The number
\[
\frac{A_k}{B_k}=\O1(q_1,q_2,\dots,q_k)=
\frac{1}{q_1}-\frac{1}{q_1q_2}+\dots
+\frac{(-1)^{k-1}}{q_1q_2\dots q_k}
\]
is called \emph{the convergent of order $k$ of the Ostrogradsky
series}.

By using the method of mathematical induction, it is easy to prove
that for any natural number $k$ the following equalities hold:
\begin{align*}
&A_k=A_{k-1}q_k+(-1)^{k-1},\\
&B_k=B_{k-1}q_k=q_1q_2\dots q_k,
\end{align*}
\textup($A_0=0, B_0=1$\textup).

From the Leibniz theorem on the convergence of alternating series,
it follows that the sequence of convergents of an even order
increases and the sequence of convergents of an odd order decreases.
Moreover, any convergent of odd order is greater than any convergent
of even order.

\section{Cylindrical sets and their properties}

\begin{definition}
A set $\bcyl1{c_1c_2\dots c_m}$, which is the closure of the set
of all  numbers $x\in(0,1)$, whose first $m$ $\bO1$-symbols are
equal to $c_1$, $c_2$, \dots, $c_m$ correspondingly, is said to be
\emph{the cylindrical set \textup(cylinder\textup) of rank $m$
with the base $(c_1, c_2, \dots, c_m)$}.
\end{definition}

Let us consider some basic properties of cylindrical sets.

\begin{property}
$\bcyl1{c_1\dots c_m}=[a,b]$, where
\begin{align*}
&a=\min\{\bO1(c_1,\dots,c_m), \bO1(c_1,\dots,c_m+1)\},\\
&b=\max\{\bO1(c_1,\dots,c_m), \bO1(c_1,\dots,c_m+1)\}.
\end{align*}
\end{property}

\textbf{Remark.} We shall denote by $\bint1{c_1\dots c_m}$ the
interior part of the set $\bcyl1{c_1\dots c_m}$.

\begin{property}
$\bcyl1{c_1\dots c_m}=\bigcup\limits_{c=1}^\infty\bcyl1{c_1\dots
c_mc}\bigcup\bO1(c_1,c_2,\dots,c_m)$, moreover
\begin{align*}
\sup\bcyl1{c_1\dots c_mc}&=\inf\bcyl1{c_1\dots c_m(c+1)}\quad
\text{if $m$ is odd},\\
\inf\bcyl1{c_1\dots c_mc}&=\sup\bcyl1{c_1\dots c_m(c+1)}\quad
\text{if $m$ is even},
\end{align*}
and
\[
\bcyl1{c_1\dots c_mc}\cap \bcyl1{c_1\dots
c_m(c+1)}=\set{\bO1(c_1,c_2,\dots,c_m,c+1)}.
\]

\end{property}

\begin{property}
$\bcyl1{c_1\dots c_m}=\bcyl1{s_1\dots s_k}$ if and only if

$m=k$ and $c_i=s_i$ for all $i=\overline{1,m}$.
\end{property}

\begin{property}
$\bcyl1{c_1\dots c_m}\subset\bcyl1{s_1\dots s_k}$ if and only if

$m\geq k$ and $c_i=s_i$ for all $i=\overline{1,k}$.
\end{property}

\begin{property}
$\bint1{c_1\dots c_m}  \cap\bint1{s_1\dots s_k}=\varnothing$ if and
only if there exists $j$ such that $c_j\not=s_j$.
\end{property}

\begin{property}\label{prop:length}
The Lebesgue measure of the cylindrical set $\bcyl1{c_1\dots c_m}$
is equal to
\[
\abs{\bcyl1{c_1\dots c_m}}=
\frac1{\sigma_1\sigma_2\dots\sigma_m(\sigma_m+1)},
\]
where $\sigma_k=\sum\limits_{i=1}^k c_i$.

\end{property}

\begin{corollary}
The cylindrical set $\bcyl1{\underbrace{\scriptstyle11\dots1}_m}$
has the largest length among the cylindrical sets of rank $m$,
namely
\[\big\lvert\bcyl1{\underbrace{\scriptstyle11\dots1}_m}\big\rvert=\frac1{(m+1)!}.\]
\end{corollary}

\textbf{Remark.} There exist cylindrical sets of different ranks
with the same lengths. For instance,
\begin{equation*}
\abs{\bcyl1{1c}}=\abs{\bcyl1 {c+1}},\quad \abs{\bcyl1{1c_2c_3\dots
c_m}}=\abs{\bcyl1{(c_2+1)c_3\dots c_m}}.
\end{equation*}

\begin{corollary}
For any given $c\in\N$ and $s\in\N $, the ratio
\[
\frac{\abs{\bcyl1{c_1\dots c_ms}}}{\abs{\bcyl1{c_1\dots c_mc}}}=
\frac{(\sigma_m+c)(\sigma_m+c+1)}{(\sigma_m+s)(\sigma_m+s+1)}
\]
converges to $1$, if $\sigma_m = \sum\limits_{i=1}^m c_i$
converges to $+\infty$.
\end{corollary}

\begin{corollary}
The ratio
\[
\frac{\abs{\bcyl1{c_1\dots c_m(c+1)}}}{\abs{\bcyl1{c_1\dots
c_mc}}}= \frac{\sigma_m+c}{\sigma_m+c+2}= 1-\frac2{\sigma_m+c+2}
\]
converges to $1$ for $m\to\infty$ \textup(or even
$\sigma_m\to\infty$\textup) or $c\to\infty$.
\end{corollary}

So, if $\sigma_m$ is large enough, then the ``weights'' of two
consecutive $\bO1$-symbols $c$ and $c+1$ are ``practically
equal''.

\section{Some metric problems and relations}

\begin{lemma}\label{lem:rel}
For any given $s \in N$, the ratio of lengths of cylindrical sets
$\bcyl1{c_1\dots c_ms}$ and $\bcyl1{c_1\dots c_m}$  satisfies the
following equality
\begin{equation}\label{eq:rel}
\frac{\abs{\bcyl1{c_1\dots c_ms}}}{\abs{\bcyl1{c_1\dots c_m}}}=
\frac{a}{(a+s-1)(a+s)}= f_s(a),
\end{equation}
where $a=1+\sigma_m$. Moreover,
\begin{equation}\label{eq:relest1}
f_s(a)\leq\frac1{2\cdot(2s-1)}
\end{equation}
and for $m\geq s-1$
\begin{equation}\label{eq:relest2}
\frac{\abs{\bcyl1{c_1\dots c_ms}}}{\abs{\bcyl1{c_1\dots c_m}}}\leq
\frac{m+1}{(m+s)(m+s+1)}.
\end{equation}
\end{lemma}

\begin{proof}
Equality~\eqref{eq:rel} follows directly from
property~\ref{prop:length} of cylindrical sets. Let us consider
\[
f_s(x)=\frac{x}{(x+s-1)(x+s)}
\]
as a function of a real variable $x$, $x\geq1$. This function
increases on $\left[1,\sqrt{(s-1)s}\right]$ and decreases on
$\left[\sqrt{(s-1)s},+\infty\right)$. Since $a$ takes only natural
values, we have
\[
\max_{a\in\N}f_s(a)=f_s(s-1)=f_s(s)=\frac1{2\cdot(2s-1)}.
\]
So, inequality \eqref{eq:relest1} holds. {The corresponding}
equality holds for $a=s$ and for $a=s-1$ (if it is possible,
because $a\geq m+1$ and it is impossible for $m\geq s$).

Function $f_s(x)$ decreases on interval $(s,+\infty)$. Hence
$f_s(a)\leq f_s(m+1)$, so inequality~\eqref{eq:relest2} holds.
\end{proof}

\textbf{Corollary.} If $c_1+\dots+c_m=s_1+\dots+s_k$ then
\[
\frac{\abs{\bcyl1{c_1\dots c_ms}}}{\abs{\bcyl1{c_1\dots c_m}}}=
\frac{\abs{\bcyl1{s_1\dots s_ks}}}{\abs{\bcyl1{s_1\dots s_k}}}.
\]

\textbf{Remark.} Let $\Delta_{c_1\dots c_m}^{\mathrm{c.f.}}$ be a
cylindrical set generated by the continued fractions representation
of real numbers. It is well known (see, e.g., \cite{Khi97}) that
\[
 \frac{\abs{\Delta_{c_1\dots
c_ms}^{\mathrm{c.f.}}}}{\abs{\Delta_{c_1\dots
c_m}^{\mathrm{c.f.}}}} = \frac{1}{s^2} \cdot \frac{1+
\frac{Q_{m-1}}{Q_m}}
 {\left(1+\frac{Q_{m-1}}{s Q_m}\right)\left(1+ \frac{1}{s} + \frac{Q_{m-1}}{s
 Q_m}\right)},
\]
where $Q_k$ is the denominator of the  $k$-th convergent  of the
continued fraction
\[
[c_1,c_1,\dots,c_n,\dots],
\]
i.e.,
\[
Q_k=c_kQ_{k-1}+Q_{k-2}\quad \text{with}\quad Q_0=1,\quad Q_1=a_1.
\]
From the latter equality it follows that the following double
inequality
\[
\frac1{3s^2}< \frac{\abs{\Delta_{c_1\dots
c_ms}^{\mathrm{c.f.}}}}{\abs{\Delta_{c_1\dots
c_m}^{\mathrm{c.f.}}}}< \frac2{s^2}
\]
holds for any sequence $(c_1,\dots,c_m)$ and for any $s\in\N$. For
the $\bO1$-representation we have $f_s(a)\to0$ ($a\to\infty$) and
Lemma~\ref{lem:rel} shows the fundamental difference between metric
relations in the representation of numbers by the first Ostrogradsky
series and by continued fractions.

\begin{lemma}\label{lem:relsum}
Let $\bcyl1{c_1\dots c_m}$ be a fixed cylindrical set, 
then
\[
\lambda\left(\bigcup_{s=1}^k \bcyl1{c_1\dots c_ms}\right)=
\frac{k}{\sigma_m+k+1}\abs{\bcyl1{c_1\dots c_m}}.
\]
\end{lemma}

\begin{proof}
From the property~\ref{prop:length} of cylindrical sets it follows
that
\begin{align*}
\lambda\left(\bigcup_{s=1}^k \bcyl1{c_1\dots c_ms}\right)
&=\sum_{s=1}^k \abs{\bcyl1{c_1\dots c_ms}}=\\
&=\frac1{\sigma_1\sigma_2\dots\sigma_m} \sum_{s=1}^k
\frac1{(\sigma_m+s)(\sigma_m+s+1)}=\\
&=\frac1{\sigma_1\sigma_2\dots\sigma_m}
\left(\frac1{\sigma_m+1}-\frac1{\sigma_m+k+1}\right)=\\
&=\frac1{\sigma_1\sigma_2\dots\sigma_m(\sigma_m+1)}\cdot\frac{k}{\sigma_m+k+1}=\\
&=\abs{\bcyl1{c_1\dots c_m}}\cdot\frac{k}{\sigma_m+k+1},
\end{align*}
which proves Lemma~\ref{lem:relsum}.
\end{proof}

\textbf{Corollary 1.} For any $k\in\N$ and for any sequence
$(c_1,\dots,c_m)$ the following inequality holds:
\[
\frac1{\sigma_m+2}\abs{\bcyl1{c_1\dots c_m}}\leq
\lambda\left(\bigcup_{s=1}^k \bcyl1{c_1\dots c_ms}\right)\leq
\frac{k}{m+k+1}\abs{\bcyl1{c_1\dots c_m}}.
\]

\textbf{Remark.} If $V\subset\N$, then it is evident that
\[
\sum_{s\in V} \abs{\bcyl1{c_1\dots c_ns}}= \abs{\bcyl1{c_1\dots
c_n}}- \sum_{s\in \N\setminus V} \abs{\bcyl1{c_1\dots c_ns}}.
\]

\textbf{Corollary 2.} Let $\bcyl1{c_1\dots c_m}$ be a fixed
cylindrical set, 
then
\[
\lambda\left(\bigcup_{c=k+1}^\infty \bcyl1{c_1\dots c_mc}\right)=
\frac{\sigma_m+1}{\sigma_m+k+1}\abs{\bcyl1{c_1\dots c_m}}.
\]

\textbf{Corollary 3.} For any $k\in\N$ and for any sequence
$(c_1,\dots,c_m)$ the following inequality holds:
\[
\frac{m+1}{m+k+1}\abs{\bcyl1{c_1\dots c_m}}\leq
\lambda\left(\bigcup_{c=k+1}^\infty \bcyl1{c_1\dots
c_mc}\right)\leq \frac{\sigma_m+1}{\sigma_m+2}\abs{\bcyl1{c_1\dots
c_m}}.
\]

\begin{theorem}
The Lebesgue measure of the set
\[
A_\sigma=\set{x : x=\bO1(g_1(x),\dots g_m(x),\dots),\;
g_{m+1}(x)>g_1(x)+\dots+g_m(x)\; \forall\, m\in\N}
\]
is equal to $0$.
\end{theorem}

\begin{proof}
Let
\[
L_k=\bigcup_{c_1\in\N}\bigcup_{c_2>\sigma_1}\dots
\bigcup_{c_k>\sigma_{k-1}} \bcyl1{c_1\dots c_k}.
\]
Then $\lambda(L_1)=\sum\limits_{c_1\in\N}\abs{\bcyl1{c_1}}=1$, and
from Corollary~2 after Lemma~\ref{lem:relsum} it follows that
\[
\lambda(L_2)= \sum_{c_1=1}^\infty \sum_{c_2=c_1+1}^\infty
\abs{\bcyl1{c_1c_2}}=\sum_{c_1=1}^\infty
\frac{c_1+1}{2c_1+1}\abs{\bcyl1{c_1}}<\frac23\sum_{c_1=1}^\infty
\abs{\bcyl1{c_1}}= \frac23\lambda(L_1),
\]
since the function $f(x)=\frac{x+1}{2x+1}$ decreases on
$(1,+\infty)$.

Similarly,
\begin{align*}
\lambda(L_{k+1}) &=\sum_{c_1=1}^\infty
\sum_{c_2=\sigma_1+1}^\infty\dots \sum_{c_{k+1}=\sigma_k+1}^\infty
\abs{\bcyl1{c_1\dots c_kc_{k+1}}}=\\
&=\sum_{c_1=1}^\infty \sum_{c_2=\sigma_1+1}^\infty\dots
\sum_{c_k=\sigma_{k-1}+1}^\infty \frac{\sigma_k+1}{2\sigma_k+1}
\abs{\bcyl1{c_1\dots c_k}}<\\
&<\frac23\sum_{c_1=1}^\infty \sum_{c_2=\sigma_1+1}^\infty\dots
\sum_{{c_k=\sigma_{k-1}+1}}^\infty \abs{\bcyl1{c_1\dots
c_k}}=\frac23\lambda(L_k).
\end{align*}
So,
\[
\lambda(L_{k+1})<\frac23\lambda(L_k)
\]
and we have
\[
\lambda(L_{k+1})<\left(\frac23\right)^k \lambda(L_1).
\]
From $A_\sigma=\bigcap\limits_{k=1}^\infty L_k$ it follows that
\[
\lambda(A_\sigma)=\lim_{k\to\infty}\lambda(L_{k+1})=0,
\]
which proves the Theorem.
\end{proof}

\textbf{Corollary.} The Lebesgue measure of set $[0,1]\setminus
A_\sigma$ is equal to $1$. That is for Lebesgue almost all  $x$:
\[g_{m+1}(x)\leq g_1(x)+\dots+g_m(x),\]
for at least  one natural $m$.

\section{The set $\Cset{V_n}$}

In this Section we shall study metric properties of the set
$\Cset{V_n}$, which is the closure of the set $\set{x: g_n(x)\in
V_n,\; n\in\N},$ consisting of the real numbers $x\in[0,1]$ whose
$\bO1$-symbols satisfy the condition
\[g_n(x)\in V_n,\]
where $\{V_n\}$ is a fixed sequence of nonempty subsets of $\N$.

It is evident that
\begin{enumerate}
\item if $V_n=\N$ for all $n\in\N$, then $\Cset{V_n}=[0,1]$,

\item if $V_n=\N$ for all $n>n_0$, then the set $\Cset{V_n}$ is a
union of segments.
\end{enumerate}

We are interested only in the case where $V_n\not=\N$ for an
infinite number of $n$.

Let $F_k = \left(\bigcup\limits_{c_1\in V_1} \ldots
\bigcup\limits_{c_k\in V_k} \bcyl1{c_1c_2\dots c_k}\right)^{\cl}$,
where $\cl$ stands for the closure.

\begin{lemma}\label{lem:CVrepresentation}
The set $\Cset{V_n}$ can be represented in the form
\[
\Cset{V_n}=\bigcap_{k=1}^\infty  F_k.
\]

It is a perfect set \textup(that is a closed set without isolated
points\textup). If $V_n\not=\N$ for an infinite number of $n$,
then it is a nowhere dense set.
\end{lemma}

\begin{proof}
The irrational number $x_0$ belongs to the set $\Cset{V_n}$ if and
only if for all natural $k$ there exists a cylindrical set
$\bcyl1{c_1c_2\dots c_k}$ of rank $k$ containing $x_0$, and $c_1\in
V_1$, $c_2\in V_2$,~\dots, $c_k\in V_k$. Let now $y_0 \in
\Cset{V_n}$ be a rational number. From the definition of the set
$\Cset{V_n}$ it follows that for any $s \in N$ the interval  $(y_0 -
\frac1s, y_0 + \frac1s)$ contains  an irrational number
$x_s\in\Cset{V_n}$ . From what has already been proved, it follows
that $x_s \in\bigcap\limits_{k=1}^\infty F_k$. Since the latter set
is closed and $x_s \rightarrow y_0$, we have $y_0 \in
\bigcap\limits_{k=1}^\infty F_k$.

The proof of the inverse inclusion is completely similar.

From
\begin{gather*}
\Cset{V_n}=[0,1]\setminus \bigcup_{m=0}^\infty
\bigcup_{\begin{subarray}{l} c_i\in V_i,\\ i=\overline{1,m}
\end{subarray}}
\Big(\bigcup_{c\in\N\setminus V_{m+1}} \bint1{c_1\dots c_mc}\cup
A_{m+1}\Big),\\
A_{m+1}=\bigcup_{i,j\in \N\setminus V_{m+1},\; i\not=j}
\left(\bcyl1{c_1\dots c_mi}\cap\bcyl1{c_1\dots c_mj}\right)
\end{gather*}
it follows that $\Cset{V_n}$ is perfect.

If $V_n\not=\N$ for an infinite number of $n$, then $\Cset{V_n}$ is
nowhere dense, because for any interval $(a,b)$ there exist an the
interval $\bint1{c_1\dots c_m}\subset(a,b)$  and  an interval
$\bint1{c_1\dots c_n\dots c_{m+k}}$ such that $ \bint1{c_1\dots
c_n\dots c_{m+k}}\cap\Cset{V_n}=\varnothing, $ where $c_{m+k}\in
\N\setminus V_{n+k}\not=\varnothing$.
\end{proof}

\textbf{Corollary.} The Lebesgue measure
\[
\lambda(\Cset{V_n})\leq \sum_{c_1\in V_1}\dots \sum_{c_k\in V_k}
\frac1{\sigma_1\dots \sigma_k(\sigma_k+1)}= \lambda(F_k)
\]
for any $k\in\N$, and
\[
\lambda(\Cset{V_n})=\lim_{k\to\infty} \lambda(F_k).
\]

Let $M_k$ be the union of all ``admissible'' cylinders of rank
$k$, i.e.,
\[
M_k = \bigcup\limits_{c_1\in V_1} \ldots \bigcup\limits_{c_k\in
V_k}~~ \bcyl1{c_1c_2\dots c_k},\quad M_0 = [0,1],
\]
and
\[
\bar M_{k+1} := M_k \setminus M_{k+1}.
\]
Then
\begin{align*}
\lambda(M_k) &=  \sum\limits_{c_1\in V_1}\ldots
\sum\limits_{c_k\in V_k} \frac1{\sigma_1\dots
\sigma_k(\sigma_k+1)},\\
\lambda (\bar M_{k+1}) &= \lambda \left( \bigcup\limits_{c_1\in
V_1} \ldots \bigcup\limits_{c_k\in V_k}~\bigcup\limits_{s\not\in
V_{k+1}} \bint1{c_1c_2\dots c_k s}\right)=\\
&= \sum\limits_{c_1\in V_1}\ldots \sum\limits_{c_k\in V_k}
\sum\limits_{s\not\in V_{k+1}} \frac1{\sigma_1\dots
\sigma_k(\sigma_k+s)(\sigma_k+s+1)} =\\
&=\sum\limits_{c_1\in V_1}\ldots \sum\limits_{c_k\in V_k} \left[
\frac1{\sigma_1\ldots \sigma_k}\sum\limits_{s\not\in
V_{k+1}}\frac1{(\sigma_k+s)(\sigma_k+s+1)} \right].
\end{align*}

\begin{lemma}\label{lem:mes0iff}
The Lebesgue measure of the set $\Cset{V_n}$ is equal to $0$ if
and only if
\[
\sum_{k=1}^\infty \frac{\lambda(\bar
M_{k+1})}{\lambda(M_k)}=+\infty.
\]
\end{lemma}

\begin{proof} Since   $F_k \setminus M_k$ is at most a countable
set, we have:
\begin{align*}
\lambda(\Cset{V_n})&=\lim_{k\to\infty}\lambda(M_{k+1})
=\lim_{k\to\infty}\frac{\lambda(M_{k+1})}{\lambda(M_k)}\cdot\frac{\lambda(M_k)}{\lambda(M_{k-1})}
\cdot ~\ldots~ \cdot\frac{\lambda(M_1)}{\lambda(M_0)} ~=\\
&=\prod_{k=0}^\infty\frac{\lambda(M_{k+1})}{\lambda(M_k)}=\prod_{k=0}^\infty\frac{\lambda(M_k)-\lambda(\bar
M_{k+1})}{\lambda(M_k)}\\
&=\prod_{k=0}^\infty\left(1-\frac{\lambda(\bar
M_{k+1})}{\lambda(M_k)}\right)=0 \iff \sum_{k=1}^\infty
\frac{\lambda(\bar M_{k+1})}{\lambda(M_k)}=+\infty,
\end{align*}
since $0<\frac{\lambda(\bar M_{k+1})}{\lambda(M_k)}<1$.
\end{proof}

First of all we shall study the problem of the determination of
the Lebesgue measure of the set $\Cset{V}=\Cset{V_n}$ with
$V_n=V$, where $V$ is a fixed proper subset of positive integers.
The sets $\Cset{V}$ with
\begin{enumerate}

\item $V=\set{1,2,\dots m}$,

\item $V=\set{m+1,m+2,\dots}$,

\item $V=\set{1,3,5,\dots}$


\end{enumerate}
are the most simple sets among $\Cset{V}$.

Let us solve the first problem in more general setting.

\begin{theorem}\label{thm:CVNk}
If the set $V_k$ contains $N_k$ symbols \textup($k\in\N$\textup)
and
\[
\varliminf_{k\to\infty}\frac{N_1N_2 \dots N_k}{(k+1)!}=0
\]
then the Lebesgue measure of the set $\Cset{V_k}$ is equal to $0$.
\end{theorem}

\begin{proof}
From the properties of  cylindrical sets it follows that
\[
\lambda(M_k)= \sum_{
\begin{subarray}{l}
v_i\in V_i,\\ i=\overline{1,k}
\end{subarray}}
\abs{\bcyl1{v_1v_2\dots v_k}}\leq \frac{N_1N_2 \dots N_k}{(k+1)!}.
\]
From Lemma~\ref{lem:CVrepresentation} and from the continuity of
Lebesgue measure it follows that
\[
\lambda(\Cset{V_k})= \lim_{k\to\infty}\lambda(M_k)\leq
\varliminf_{k\to\infty}\frac{N_1N_2 \dots N_k}{(k+1)!}=0.\qed
\]
\renewcommand{\qed}{\relax}
\end{proof}

\textbf{Corollary.} If $N_k\leq m$ \textup(for any
$k\in\N$\textup) for some fixed $m$, then the Lebesgue measure of
the set $\Cset{V_k}$ is equal to $0$.

\begin{theorem} \label{theorem4}
Let $V_k=\set{1,2,\dots,m_k}$. If $\sum\limits_{k=1}^\infty
\frac1{m_k}=+\infty$, then the Lebesgue measure of the set
$\Cset{V_n}$ is equal to $0$.
\end{theorem}

\begin{proof}
Let $\bcyl1{c_1c_2\dots c_k}$ be a fixed cylindrical set of rank
$k$. Then
\begin{align*}
\sum_{c\not\in V_{k+1}} \abs{\bint1{c_1c_2\dots
c_kc}}&=\frac1{\sigma_1\sigma_2\dots\sigma_k}
\sum_{c=m_k+1}^\infty
\frac1{(\sigma_k+c)(\sigma_k+c+1)}=\\
&=\frac1{\sigma_1\sigma_2\dots\sigma_k(\sigma_k+m_k+1)}.
\end{align*}
Since
\[
\frac1{\sigma_k+m_k+1}>\frac1{(m_k+1)(\sigma_k+1)},
\]
we have
\[
\sum_{c\not\in V_{k+1}} \abs{\bint1{c_1c_2\dots c_kc}}>
\frac1{m_k+1}\cdot \abs{\bcyl1{c_1c_2\dots c_k}}.
\]
Summing over all $c_1\in V_1$, $c_2\in V_2$, \dots, $c_k\in V_k$,
we have
\[
\lambda(\bar M_{k+1})>\frac1{m_k+1}\lambda(M_k),\quad
\text{i.e.},\quad \frac{\lambda(\bar
M_{k+1})}{\lambda(M_k)}>\frac1{m_k+1}
\]
for any $k\in\N$, and the statement of the Theorem follows
directly from Lemma~\ref{lem:mes0iff}.
\end{proof}

Let $E$ be the set of all real numbers with bounded
$\bO1$-symbols, i.e.,  $x\in E$ iff there exists a constant $K_x$
such that $g_k(x)\leq K_x$ for all $k\in N$.

\begin{theorem}
The Lebesgue measure of the set $E$ of all real numbers $x\in[0,1]$
with bounded $\bO1$-symbols is equal to $0$.
\end{theorem}

\begin{proof}
For a given $m\in N$, let us consider the set $E_m=\{x: g_k(x)\leq
m, \forall k\in N\}$ of $m$-uniformly bounded symbols. It is not
hard to see that $E_m = \Cset{V_k}$ with $V_k = \{1,2, ..., m \}$.
From the latter Theorem it follows that $\lambda(\Cset{V_k})=0$.

Since  $E=\bigcup\limits_{m=1}^{\infty} E_m$ and $\lambda(E_m)=0$,
we have the desired conclusion.
\end{proof}

\textbf{Corollary.} For Lebesgue almost all real numbers $x\in[0,1]$
the following equality holds:
\[
\varlimsup_{k\to\infty} g_{k}(x)=\infty.
\]


Let us now consider the case, where $V_k = \set{v_k+1,v_k+2,\dots},$
and $\{v_k \}$ is a fixed sequence of positive integers.

\begin{lemma}\label{lemma5}
Let $\bcyl1{c_1\dots c_n}$ be a fixed cylindrical set or, if $n=0$,
the unit interval $[0,1]$;

 let $\{v_k \}$ be a fixed sequence of
positive integers, let $V_k=\set{v_k+1,v_k+2,\dots}$, and let
\begin{align*}
 M^{c_1 \dots c_n}_k &:= M_{n+k} \bigcap
\bcyl1{c_1\dots c_n} =\bigcup_{c_{n+1}>v_{n+1}} \dots
\bigcup_{c_{n+k}>v_{n+k}} \bcyl1{c_1\dots
c_nc_{n+1}\dots c_{n+k}},\\
\bar M^{c_1 \dots ~c_n}_{k+1}&=M^{c_1 \dots c_n}_k\setminus M^{c_1
\dots c_n}_{k+1}=\bigcup_{c_{n+1}>v_{n+1}} \dots
\bigcup_{c_{n+k}>v_{n+k}} \bigcup_{s=1}^{v_{n+k+1}}
\bint1{c_1\dots c_{n+k}s}.
\end{align*}
Then
\begin{equation}
\frac{\lambda(\bar M^{c_1 \dots c_n}_{k+1})}{\lambda(\bar M^{c_1
\dots c_n}_k)} <\frac12\cdot\frac{v_{n+k+1}}{v_{n+k}}.
\end{equation}
\end{lemma}

\begin{proof}
Let $\bint1{c_1\dots c_n\dots c_{n+k-1}s}$ be a fixed cylindrical
interval of rank $n+k$. Then
\begin{align*}
\sum_{s\not\in V_{n+k}}\abs{\bint1{c_1\dots c_n\dots c_{n+k-1}s}}&
= \sum\limits_{s=1}^{v_{n+k}}
\frac1{\sigma_1\dots\sigma_{n+1}\ldots\sigma_{n+k-1}(\sigma_{n+k-1}+s)(\sigma_{n+k-1}+s+1)}=\\
&=\frac1{\sigma_1\dots\sigma_{n+1}\dots \sigma_{n+k-1}}
\left(\frac1{\sigma_{n+k-1}+1}-\frac1{\sigma_{n+k-1}+
v_{n+k}+1}\right)=\\
&=\frac{v_{n+k}}{\sigma_1\dots\sigma_{n+1}\dots \sigma_{n+k-1}
(\sigma_{n+k-1} +1)(\sigma_{n+k-1} +v_{n+k}+1)}.
\end{align*}
Let $\bint1{c_1\dots c_n\dots c_{n+k}s}$ be a fixed cylindrical
interval of rank $n+k+1$. Then
\begin{align*}
\sum_{c_{n+k}\in V_{n+k}}& \sum_{s\not\in V_{n+k+1}}
\abs{\bint1{c_1\dots c_n\dots c_{n+k}s}} =\\
&=\sum_{c_{n+k}=v_{n+k}+1}^\infty
\frac1{\sigma_1\dots\sigma_{n+1}\dots \sigma_{n+k}}
\left(\frac1{\sigma_{n+k}+1}-\frac1{\sigma_{n+k}+
v_{n+k+1}+1}\right)=\\
&=\frac1{\sigma_1\dots\sigma_{n+1}\dots \sigma_{n+k-1}}
\sum_{s=v_{n+k}+1}^\infty
\biggl(\frac1{(\sigma_{n+k-1}+s)(\sigma_{n+k-1}+s+1)}-{}\\
&\qquad\qquad\qquad\qquad\qquad-\frac1{(\sigma_{n+k-1}+s)(\sigma_{n+k-1}+s+v_{n+k+1}+1)}\biggr) = \\
&=\frac1{\sigma_1\dots\sigma_{n+1}\dots\sigma_{n+k-1}}
\biggl(\frac1{\sigma_{n+k-1}+v_{n+k}+1} - {} \\
&\qquad\qquad\qquad\qquad\qquad-\frac1{1+v_{n+k+1}}
\sum_{i=1}^{v_{n+k+1}+1} \frac1{\sigma_{n+k-1}+v_{n+k}+ i}
\biggr) =\\
&=\frac{v_{n+k}}{\sigma_1\dots\sigma_{n+1}\dots\sigma_{n+k-1}(\sigma_{n+k-1}+1)(\sigma_{n+k-1}+v_{n+k}+1)}\cdot
X_k.
\end{align*}

Let us estimate the expression
\begin{align*}
X_k
&=\frac{(\sigma_{n+k-1}+1)(\sigma_{n+k-1}+v_{n+k}+1)}{v_{n+k}}\cdot
\Biggl( \frac1{\sigma_{n+k-1}+v_{n+k}+1}-{} \\
&\qquad\qquad\qquad\qquad\qquad-\frac1{1+v_{n+k+1}}\cdot
\sum_{i=1}^{1+v_{n+k+1}}\frac1{\sigma_{n+k-1}+v_{n+k}+i}\Biggr) =\\
&=\frac{\sigma_{n+k-1}+1}{v_{n+k}}\cdot
\left(1-\frac1{1+v_{n+k+1}}\cdot
\sum_{i=1}^{1+v_{n+k+1}}\frac{\sigma_{n+k-1}+v_{n+k}+1}{\sigma_{n+k-1}+v_{n+k}+i}\right) =\\
&=\frac{\sigma_{n+k-1}+1}{v_{n+k}}\cdot
\left(1-\frac1{1+v_{n+k+1}}\cdot
\sum_{i=1}^{1+v_{n+k+1}}\left(1-\frac{i-1}{\sigma_{n+k-1}+v_{n+k}+i}\right)\right) = \\
&=\frac{\sigma_{n+k-1}+1}{v_{n+k}}\cdot \frac1{1+v_{n+k+1}}\cdot
\sum_{i=2}^{1+v_{n+k+1}}\frac{i-1}{\sigma_{n+k-1}+v_{n+k}+i}.
\end{align*}
Now let us estimate the following sum
\[
\frac1{n_0+1}+\frac2{n_0+2}+\dots+\frac{m_k}{n_0+m_k},
\]
where $n_0$ and $m_k>1$~ are natural numbers.  Let
\[
C_k := \frac1{n_0+1}+\frac1{n_0+2}+\dots+\frac1{n_0+m_k}
\]
and let us consider the following ``matrix'':
\[
\begin{array}{ccccc}
\frac1{n_0+1}&\frac1{n_0+2}&\frac1{n_0+3}&\dots&\frac1{n_0+m_k}\\
\frac1{n_0+1}&\frac1{n_0+2}&\frac1{n_0+3}&\dots&\frac1{n_0+m_k}\\
\frac1{n_0+1}&\frac1{n_0+2}&\frac1{n_0+3}&\dots&\frac1{n_0+m_k}\\
\vdots&\vdots&\vdots&\vdots&\vdots\\
\frac1{n_0+1}&\frac1{n_0+2}&\frac1{n_0+3}&\dots&\frac1{n_0+m_k}\\
\end{array}
\]

The sum of all addends over the whole ``matrix'' is equal to
$m_k\cdot C_k$.

The sum of all addends over the ``main diagonal of the matrix'' is
equal to $C_k$.

 The sum of all elements standing above ``the main diagonal''
 is less then the sum of all elements standing under ``the main diagonal''
 (for any element above the ``main diagonal'' there exists the symmetrical element (under the ``main diagonal''),
 which is greater then the initial one).

  The sum of all elements standing outside  the ``main diagonal'' is equal to $(m_k-1)\cdot C_k$.
So, the sum of all elements standing above the ``main diagonal''
is less than $\frac{m_k-1}2\cdot C_k$, and the sum of all elements
above the ``main diagonal'' and over the ``main diagonal''  is
equal to
\[
\frac1{n_0+1}+\frac2{n_0+2}+\dots+\frac{m_k}{n_0+m_k}<\frac{m_k-1}2\cdot
C_k+C_k=\frac{m_k+1}2\cdot C_k.
\]
So,
\[
\frac1{n_0+1}+\frac2{n_0+2}+\dots+\frac{m_k}{n_0+m_k}<\frac{m_k+1}2\cdot
\left(\frac1{n_0+1}+\frac1{n_0+2}+\dots+\frac1{n_0+m_k}\right).
\]
Therefore,
\begin{align*}
X_k&=\frac{\sigma_{n+k-1}+1}{v_{n+k}}\cdot
\frac1{1+v_{n+k+1}}\cdot
\sum_{i=1}^{v_{n+k+1}}\frac{i}{(\sigma_{n+k-1}+v_{n+k}+1)+i}<\\
&<\frac{\sigma_{n+k-1}+1}{v_{n+k}}\cdot \frac1{1+v_{n+k+1}}\cdot
\frac{v_{n+k+1}+1}2\cdot \sum_{i=1}^{v_{n+k+1}}\frac1{\sigma_{n+k-1}+v_{n+k}+i+1}=\\
&=\frac1{2v_{n+k}}
\sum_{i=1}^{v_{n+k+1}}\frac{\sigma_{n+k-1}+1}{\sigma_{n+k-1}+v_{n+k}+i+1}<\frac12\cdot\frac{v_{n+k+1}}{v_{n+k}}.
\end{align*}

So, the inequality
\[
\sum_{c_{n+k}\in V_{n+k}} \sum_{s\not\in V_{n+k+1}}
\abs{\bint1{c_1\dots c_n\dots c_{n+k}s}} <
\frac12\cdot\frac{v_{n+k+1}}{v_{n+k}}\cdot \sum_{s\not\in
V_{n+k}}\abs{\bint1{c_1\dots c_n\dots c_{n+k-1}s}}
\]
holds. Hence, summing over all $c_{n+1}\in V_{n+1}$, $c_{n+2}\in
V_{n+2}$, $\dots$, $c_{n+k-1}\in V_{n+k-1}$, we have
\[
\lambda(\bar M^{c_1 \dots c_n}_{k+1})
<\frac12\cdot\frac{v_{n+k+1}}{v_{n+k}}\cdot\lambda(\bar M^{c_1
\dots c_n}_k),
\]
which proves the Lemma.
\end{proof}

\textbf{Corollary 1.} Let $V_k=\set{v_k+1,v_k+2,\dots}$, $v_k\in\N$.
Then
\[
\lambda(\bar M_{k+1}) <\frac12\cdot\frac{v_{k+1}}{v_k} \lambda(\bar
M_k).
\]

\textbf{Corollary 2.} Let $V_k = V=\set{m+1, m+2,\dots}$,
$m\in\N$. Then
\begin{equation*}
\lambda\left(\bar M_{k+1}^{c_1 c_2\dots c_n}\right)
<\frac12\lambda\left(\bar M_k^{c_1 c_2\dots c_n}\right)
\end{equation*}
for any natural number $k$ and any $c_1 \in V$, $\dots$, $c_n \in
V$, and, therefore,
\[
\lambda(\bar M_{k+1}) <\frac12 \lambda(\bar M_k).
\]

\begin{theorem}\label{c_0<2}
Let $\{v_k\}$ be a fixed sequence of positive integers, and let
\[
V_k=\set{v_k+1,v_k+2,\dots}.
\]

If there exists $k_0\in\N$ such that
\[
\frac{v_{k+1}}{v_k}\leq C_0 < 2 \quad \text{for any $k > k_0$},
\]
then the set $\Cset{V_k}$ is of positive Lebesgue measure.
\end{theorem}

\begin{proof}
Let $\bcyl1{c_1\dots c_n}$ be any fixed cylindrical set with
$n>k_0$ and $c_i \in V_i$. We shall prove that the set
\[
\Delta_{c_1\dots c_n}=\Cset{V_k}\cap\bcyl1{c_1\dots c_n}
\]
has positive Lebesgue measure. To this aim let us consider a
cylindrical set $\bcyl1{c_1\dots c_n c_{n+1}}$, $c_{n+1}>
v_{n+1}$, and the corresponding subset
\[
\Delta_{c_1\dots c_n c_{n+1}}=\Cset{V_k}\cap \bcyl1{c_1\dots c_n
c_{n+1}}.
\]
From Lemma~\ref{lemma5} it follows that
\begin{align*}
\lambda&(\bar M_{k+1}^{c_1 \dots c_n c_{n+1}}) <
\frac12\cdot\frac{v_{n+k+1}}{v_{n+k}}\cdot\lambda(\bar M_k^{c_1
\dots c_n c_{n+1}}) \leq \frac12\cdot C_0 \cdot\lambda(\bar
M_k^{c_1
\dots c_n c_{n+1}}) <\\
&< \frac12\cdot C_0 \cdot
\frac12\cdot\frac{v_{n+k}}{v_{n+k-1}}\cdot\lambda(\bar
M_{k-1}^{c_1 \dots c_n c_{n+1}}) \leq \left(\frac{C_0}{2}\right)^2
\cdot\lambda(\bar M_{k-1}^{c_1
\dots c_n c_{n+1}}) < \dots \\
&\leq \left(\frac{C_0}{2}\right)^k \cdot\lambda(\bar M_1^{c_1
\dots c_n c_{n+1}})
\end{align*}
for any $k\in\N$. Using Lemma~\ref{lem:relsum}, we have
\[
\lambda(\bar M^{c_1 \dots c_n c_{n+1}}_1) =
\sum\limits_{s=1}^{v_{n+2}} \abs{\bint1{c_1\dots c_n c_{n+1}s}} =
\frac{v_{n+2}}{\sigma_{n+1} + v_{n+2}+1} \cdot
\abs{\bcyl1{c_1\dots c_n c_{n+1}}}.
\]
So,
\begin{align*}
\lambda(\Delta_{c_1\dots c_n c_{n+1}})&= \abs{\bcyl1{c_1\dots c_n
c_{n+1}}} - \sum\limits_{k=1}^{\infty}
\lambda(\bar M^{c_1 \dots c_n c_{n+1}}_k)>\\
& > \abs{\bcyl1{c_1\dots c_n c_{n+1}}} -
\sum\limits_{k=1}^{\infty} \left(\frac{C_0}{2}\right)^{k-1}
\cdot\lambda(\bar M^{c_1 \dots c_n c_{n+1}}_{1}) = \\
& = \abs{\bcyl1{c_1\dots c_n c_{n+1}}} \cdot \left( 1 -
\frac{2}{2-C_0}\cdot \frac{v_{n+2}}{\sigma_{n+1}+
v_{n+2}+1}\right).
\end{align*}
Since the numbers $c_1$, \dots, $c_n$, $v_{n+2}$, $C_0$ are fixed,
and $c_{n+1}> v_{n+1}$, there exists a number $c^* \in\N$ such
that
\[
1 - \frac{2}{2-C_0}\cdot \frac{v_{n+2}}{\sigma_{n+1}+ v_{n+2}+1} >
0
\]
for any $c_{n+1}> c^*$. Hence, $\lambda(\Delta_{c_1\dots c_n
c_{n+1}})>0$ for any $c_{n+1}> c^*$, and, therefore,
\[
\lambda(\Cset{V_k}) > \lambda(\Delta_{c_1\dots c_n}) >
\lambda(\Delta_{c_1\dots c_n c_{n+1}}) >0.\qedhere
\]
\end{proof}

\textbf{Corollary.} Let $P_n(x)= a_n x^n + a_{n-1} x^{n-1} + \dots
+a_1 x^1 + a_0$ with $n\in\N$, $a_i\in\Z$ and $P_n(x)>0$ for any
$x\in\N$. If $v_k = P_n(k),$ then $\lambda(\Cset{V_k})>0$.

\begin{theorem}\label{thm:Csetnot1..m}
Let $m$ be a fixed natural number and
$V=\N\setminus\set{1,2,\dots,m}$, then the set $\Cset{V}$ is of
positive Lebesgue measure and
\begin{equation}\label{eq:Csetnot1..m.mes.est}
\lambda(\Cset{V})>\frac1{(m+1)^2}.
\end{equation}
\end{theorem}

\begin{proof}
The first statement of the Theorem follows directly from the
Theorem~\ref{c_0<2}. Let us prove the second statement. To this
aim we consider an arbitrary cylindrical set $\bcyl1{c_1c_2\dots
c_m}$ such that $c_1\in V$, $c_2\in V$,~\dots, $c_m\in V$. From
the Corollary 2 after Lemma \ref{lemma5} it follows that
\[
\lambda\left(\bar
M_{k+1}^{c_1c_2\dots c_m}\right) <\frac1{2^k}\lambda\left(\bar
M_1^{c_1c_2\dots c_m}\right).
\]
So, we have
\begin{align*}
\lambda(\Delta_{c_1c_2\dots c_m})&=\abs{\bcyl1{c_1c_2\dots
c_m}}-\sum\limits_{k=1}^\infty \lambda\left(\bar M_k^{c_1c_2\dots
c_m}\right)> \\
&>\abs{\bcyl1{c_1c_2\dots c_m}}-\lambda\left(\bar M_1^{c_1c_2\dots
c_m} \right)\cdot \sum\limits_{k=0}^\infty \frac1{2^k} =
\abs{\bcyl1{c_1c_2\dots c_m}}-2\lambda\left(\bar M_1^{c_1c_2\dots
c_m}\right).
\end{align*}
Since
\[
\lambda\left(\bar M_1^{c_1c_2\dots c_m}\right) =
\sum\limits_{c=1}^m \bint1{c_1c_2\dots c_mc} =
\frac{m}{\sigma_m+m+1}\cdot\abs{\bcyl1{c_1c_2\dots c_m}}
\leq\frac{m}{(m+1)^2}\cdot\abs{\bcyl1{c_1c_2\dots c_m}},
\]
it follows that
\[
\lambda(\Delta_{c_1c_2\dots
c_m})>\frac{m^2+1}{(m+1)^2}\cdot\abs{\bcyl1{c_1c_2\dots c_m}}.
\]

Now we shall estimate the Lebesgue measure of
$\bigcup\limits_{c_1\in V}\dots\bigcup\limits_{c_m\in
V}\abs{\bcyl1{c_1c_2\dots c_m}}:$
\begin{align*}
\sum_{c_1=m+1}^\infty\dots&\sum_{c_m=m+1}^\infty
\frac1{\sigma_1\sigma_2\dots\sigma_{m-1}(\sigma_{m-1}+c_m)(\sigma_{m-1}+c_m+1)}=\\
&=\sum_{c_1=m+1}^\infty\dots\sum_{c_{m-1}=m+1}^\infty
\frac1{\sigma_1\sigma_2\dots\sigma_{m-1}(\sigma_{m-1}+m+1)}>\\
&>\sum_{c_1=m+1}^\infty\dots\sum_{c_{m-1}=m+1}^\infty
\frac1{\sigma_1\sigma_2\dots\sigma_{m-2}(\sigma_{m-1}+m)(\sigma_{m-1}+m+1)}=\\
&=\dots=\sum_{c_1=m+1}^\infty
\frac1{\sigma_1(\sigma_1+(m-1)m+1)}>\\
&>\sum_{c_1=m+1}^\infty
\frac1{(c_1+(m-1)m)(c_1+(m-1)m+1)}=\frac1{m^2+1}.
\end{align*}
Since
\[
\lambda(\Cset{V})=\sum_{c_1\in V}\dots\sum_{c_m\in V}
\lambda(\Delta_{c_1c_2\dots c_m}),
\]
we have inequality~\eqref{eq:Csetnot1..m.mes.est}.
\end{proof}

\textbf{Corollary.} Let the sequence $\{v_k\}$ be uniformly
bounded (i.e., there exists a number $D_0 $ such that $v_k \leq
D_0, \forall~ k \in N)$. Then the set $\Cset{V_k}$ is of positive
Lebesgue measure.

\bigskip

Finally, let us consider the more general case where
$V_k=V=\N\setminus\set{a_1,a_2,\dots,a_n,\dots}$ and $\set{a_n}$ is
an arbitrary increasing sequence of positive integers.

\begin{theorem}\label{theorem: <d}
Let $\{a_n\}$ be an increasing sequence of positive integers with
$a_{n+1}-a_n\leq d$ for some fixed natural number $d\geq2$, and
for any $n\in\N$. If~
$V_k=V=\N\setminus\set{a_1,a_2,\dots,a_n,\dots}$,  then the
Lebesgue measure of the set $\Cset{V}$ is equal to $0$.
\end{theorem}

\begin{proof}
Let us fix a cylindrical set $\bcyl1{c_1c_2\dots c_k}$ and
estimate the following  sum
\begin{align*}
\sum_{c\not\in V} \abs{\bint1{c_1c_2\dots c_kc}}&=
\frac1{\sigma_1\sigma_2\dots\sigma_k} \sum_{n=1}^\infty
\frac1{(\sigma_k+a_n)(\sigma_k+a_n+1)}>\\
&>\frac1{\sigma_1\sigma_2\dots\sigma_k} \sum_{n=1}^\infty
\frac1{(\sigma_k+a_n')(\sigma_k+a_n'+d)}=\\
&=\frac1d\cdot
\frac1{\sigma_1\sigma_2\dots\sigma_k(\sigma_k+a_1)},
\end{align*}
where $a_1'=a_1$, $a_{n+1}'=a_n'+d\geq a_{n+1}$ for all natural
$n$. Since
\[
\frac1{\sigma_k+a_1}\geq\frac1{a_1(\sigma_k+1)},
\]
we have
\[
\sum_{c\not\in V} \abs{\bint1{c_1c_2\dots c_kc}}>
\frac1{a_1d}\cdot \abs{\bcyl1{c_1c_2\dots c_k}}.
\]
Summing over all $c_1\in V$, $c_2\in V$, \dots, $c_k\in V$, we
have
\[
\lambda(\bar M_{k+1})>\frac1{a_1d}\lambda(M_k),\quad
\text{i.e.},\quad \frac{\lambda(\bar
M_{k+1})}{\lambda(M_k)}>\frac1{a_1d}
\]
for any $k\in\N$, and the statement of the Theorem follows
directly from Lemma~\ref{lem:mes0iff}.
\end{proof}

\textbf{Corollary 1.}\label{lem:CVsequence} If
$V_k=V=\set{b_1,b_2,\ldots,b_n,\ldots}$ with $ b_{n+1}-b_n\geq2$,
then the Lebesgue measure of the set $\Cset{V}$ is equal to $0$.

\textbf{Corollary 2.} If $V=\set{1,3,5,\dots}$ or
$V=\set{2,4,6,\dots}$ then $\lambda(\Cset{V})=0$.

\section{Random variables with independent $\bO1$-symbols}

Let us consider the following random variable
\[
\xi=\bO1(\xi_1,\xi_2,\dots,\xi_k,\dots) :=
\sum\limits_{n=1}^{\infty}\frac{(-1)^{n-1}}{\xi_1(\xi_1+\xi_2)\cdot
\ldots \cdot(\xi_1+\xi_2+\dots+\xi_n)},
\]
where $\xi_k$ are independent random variables taking the values
$1$, $2$,~$\dots$, $m$,~$\dots$ with probabilities $p_{1k}$,
$p_{2k}$,~$\dots$, $p_{mk}$,~$\dots$ correspondingly, $p_{mk}\geq
0$, $\sum\limits_{m=1}^\infty p_{mk}=1$.

Since  the random variable $\xi$ is a sum of an infinite number of
terms, it can takes irrational values only.

\begin{lemma}[\cite{PB04}]
The distribution function $F_\xi$ of the random variable $\xi$ is
of the following form
\begin{equation}\label{eq:distributionfunctionxi}
F_\xi(x)=\beta_1(x)+\sum_{k\geq2}\left((-1)^{k-1}\beta_k(x)\prod_{i=1}^{k-1}p_{g_i(x)i}\right),\quad
\text{if}\quad 0<x\leq1,
\end{equation}
where
\[
\beta_k(x)=1-\sum_{j=1}^{g_k(x)-1}p_{jk},
\]
and $g_k(x)$ is the $k$-th $\bO1$-symbol of number $x$, and the
expression~\eqref{eq:distributionfunctionxi} has a finite resp.
infinite number of terms according to rationality resp.
irrationality of the number $x$.
\end{lemma}

\begin{theorem}\label{thm:discretenesscriterion}
The random variable $\xi$ has a discrete distribution if and only
if
\begin{equation}\label{eq:discretenesscriterion}
\prod_{k=1}^\infty \max_m\left\{p_{mk}\right\}>0,
\end{equation}
and it is  continuously distributed if and only if the infinite
product in \eqref{eq:discretenesscriterion} diverges to $0$.

Moreover, in the discrete  case the atomic spectrum of the
distribution of the random variable $\xi$ consists of real numbers
$x\in[0,1]$ whose $\bO1$-representation differs from the
$\bO1$-representation of
\[
x_0=\bO1(g'_1,g'_2,\dots,g'_k,\dots)\quad \text{with}\quad p_{g'_k
k}=\max_m\left\{p_{mk}\right\} \quad \text{for all}\quad k \in\N,
\]
by at most a finite number of $\bO1$-symbols $g_k(x)$ with
$p_{g_k(x)k}>0$.
\end{theorem}

\begin{proof}
If $\xi$ has an atomic distribution, then there exists a point $x$
such that
\[
P \{\xi = x\}=\prod_{k=1}^\infty p_{g_k(x)k}>0.
\]
 In such a
case we have
$$\prod_{k=1}^\infty \max_m\left\{p_{mk}\right\} \geq
\prod_{k=1}^\infty p_{g_k(x)k}> 0.$$
 So, if  the distribution of
the random variable $\xi$ has atoms, then
\eqref{eq:discretenesscriterion} holds.

Let now \eqref{eq:discretenesscriterion} holds. Let us consider an
arbitrary $x\in[0,1]$ whose $\bO1$-representation differs from the
$\bO1$-representation of  the above $x_0$ by at most a finite number
of $\bO1$-symbols $g_k(x)$ with $p_{g_k(x)k}>0$. It is evident that
$x$ is also an atom of the distribution $\xi$. We shall prove that
$\xi$ has a discrete distribution.

Let $x_j^{(m)}=\bO1(g_1,g_2,\dots, g_m,g'_{m+1},\dots,
g'_k,\dots)$ be an  arbitrary atom among all atoms whose
$\bO1$-symbols coincide with the $\bO1$-symbols of $x_0$ starting
from the  $(m+1)$-th symbol. Then
\begin{align*}
P\set{\xi\in\set{x_j^{(m)}}}&= \sum_{\begin{subarray}{l}
 g_1:p_{g_11}>0\\
\dots\dots\dots\dots\dots\dots\\
g_m :p_{g_mm}>0
\end{subarray}}
{\left(p_{g_11}p_{g_22} \dots
p_{g_mm}\prod_{k=m+1}^\infty{p_{g'_k(x)k}}\right)} =\\
&=\sum_{g_1:p_{g_11}>0}{p_{g_11}} \sum_{g_2:p_{g_22}>0}{p_{g_22}}
\dots \sum_{g_m:p_{g_mm}>0}{p_{g_mm}}
\prod_{k=m+1}^\infty{p_{g'_k(x)k}}=\\
&= \prod_{k=m+1}^\infty{p_{g'_k(x)k}}.
\end{align*}

The set $D=\bigcup\limits_{m=1}^\infty \set{x_j^{(m)}}$ is at  most
a countable set and
\[
P\set{\xi \in D}=\lim_{m\to\infty} P\set{\xi\in \set{x_j^{(m)}}}=
\lim_{m\to\infty} \prod_{k=m+1}^\infty{p_{g'_kk}}=1.
\]
So, the  random variable $\xi$ is supported by an at most
countable set and thus it is  discretely distributed by
definition.
\end{proof}

\begin{theorem}
The distribution of the random variable $\xi$ is of pure type. It
is either  pure discrete or pure singular continuous or pure
absolutely continuous.
\end{theorem}

\begin{proof}
Taking  into account  Theorem~\ref{thm:discretenesscriterion}, it
is sufficient to prove that in the continuous case the
distribution of $\xi$ is either  pure singular or pure absolutely
continuous.

Let $x=\bO1(g_1(x),g_2(x),\dots,g_n(x),\dots)$ and let $t_1$, \dots,
$t_n$ be fixed natural numbers. We shall set
\[
\bar\Delta_{t_1\dots t_n}(x)=
\bO1(t_1,\dots,t_n,g_{n+1}(x),g_{n+2}(x),\dots)
\]
and for any set $E\subset[0,1]$ we shall set
\begin{gather*}
\bar\Delta_{t_1\dots t_n}(E)=\set{u :
u=\bar\Delta_{t_1\dots t_n}(x), x\in E},\\
T_n(E)=\bigcup_{t_1,\dots,t_n} \bar\Delta_{t_1\dots t_n}(E),\quad
T(E)=\bigcup_n T_n(E).
\end{gather*}

Let us consider an event $A=\set{\xi\in T(E)}$. Since the random
variables $\xi_k$ are independent, the event $A$ is residual.
 So, from the Kolmogorov $0$--$1$ law it follows that either $P(A)=0$ or $P(A)=1$.

Since $T(E) \supset E$, from the inequality $P\set{\xi \in E}>0$ it
follows that $P\set{\xi \in T(E)} \geq P\set{\xi \in E}>0$, so
$P\set{\xi \in T(E)}=1$.

Only one of the following two cases can occur:
\begin{enumerate}

\item There exists a set $E$ such that $\lambda(E)=0$, but $P\set{\xi
\in E}>0$.

\item For any set $E$ with $\lambda(E)=0$ it follows that $P\set{\xi
\in E}=0$.

\end{enumerate}

In the first case from equality $\lambda(E)=0$ it follows that
$\lambda(T(E))=0$, which implies that there exists a set $T(E)$ such
that $\lambda(T(E))=0$, but $P\set{\xi \in T(E)}=1$, that is the
distribution of $\xi$ is pure singular by definition.

In the second case the distribution  of the random variable $\xi$ is
absolutely continuous by definition.
\end{proof}

Now let us consider metric and topological properties of the
topological support (i.e., the minimal closed support) $S_{\xi}$
of the random variable $\xi$. These properties are  completely
determined by the infinite stochastic matrix $P_{\xi} =
\|p_{ik}\|$, where the $k$-th column of the matrix corresponds to
the distribution of the random variable $\xi_k$: $p_{ik}=
P\{\xi_k=i\}$.

\begin{theorem} \label{Contorovity}
The topological support  $S_{\xi}$ of the random variable $\xi$ is a
nowhere dense set if and only if the matrix $P_{\xi}$ contains an
infinite number of columns having zero elements.

If the set $V_k (\xi)= \set{i: p_{ik}>0}$ has one of the following
properties:
\begin{enumerate}
\item $V_k (\xi)$ contains $N_k$ elements and $\varliminf\limits
_{k\to\infty}\frac{N_1N_2 \dots N_k}{(k+1)!}=0$;

\item $V_k(\xi)=\set{1,2,\dots,m_k}$, and
$\sum\limits_{k=1}^\infty \frac1{m_k}=+\infty$;

\item $V_k(\xi)=V=\N\setminus\set{a_1,a_2,\dots,a_n,\dots}$, where
$a_n$ is an arbitrary increasing sequence of positive integers
with $a_{n+1}-a_n\leq d$ for some fixed  $d\geq2$ and for any
$n\in\N$;
\end{enumerate}
then the topological support of the random variable $\xi$ is of
zero Lebesgue measure.

\end{theorem}

\begin{proof}
It is well known that for any arbitrary random variable $\eta$
with the distribution function $F_{\eta}$ the topological support
$S_{\eta}$ coincides with the set
\[
\set{x:
F_{\eta}(x+\varepsilon)- F_{\eta}(x-\varepsilon) > 0, \forall
\varepsilon>0}.
\]

Let us consider the set $\Cset{V_k(\xi)}$  with $V_k(\xi)=\set{i:
p_{ik}>0}$.

If $x= \bO1(g_1(x),g_2(x),\dots,g_n(x),\dots) \in \Cset{V_k(\xi)}$,
then
\[
P\{\xi \in \bcyl1{g_1(x)g_2(x)\dots g_n(x)} \}= \prod\limits_{k=1}^n
p_{g_k(x)k}>0,
\]
for any $n\in\N$. So, $x \in S_{\xi}$.

If $x= \bO1(g_1(x),g_2(x),\dots,g_n(x),\dots) \not\in
\Cset{V_k(\xi)}$, then there exists a number $n_0$ such that
$g_{n_0}(x) \not\in V_{n_0}(\xi).$ So, $p_{g_{n_0}(x) n_0}=0$, and
\[
P\{\xi \in \bcyl1{g_1(x)g_2(x)\dots g_{n_0}(x)} \}=
\prod\limits_{k=1}^{n_0} p_{g_k(x)k} =0.
\]
Hence, $x \not\in S_{\xi}$. Therefore, the topological support
$S_\xi$ of the random variable $\xi$ coincides with the set
$\Cset{V_k(\xi)}$.

If the  matrix $P_{\xi}$ contains only a finite number  of columns
having zero elements (i.e., there exists a number $k_0$ such that
$p_{ik}>0$ for any $k>k_0$ and for any $i\in\N$), then the
topological support $S_{\xi}$ completely contains any cylindrical
set $\bcyl1{c_1\dots c_k} $ with $k>k_0$ and $c_i \in V_i.$

If the  matrix $P_{\xi}$ contains an infinite number  of columns
having zero elements, then for any $n\in \N$ there exists a column
$l_n>n$ and a number $s_n \in \N$, such that $p_{s_n l_n}=0$.
Therefore, for any cylindrical set $\bcyl1{c_1c_2\dots c_n}$ with
$c_i \in V_i $ there exists a subset $\bcyl1{c_1c_2\dots c_n \dots
c_{l_n-1} s_n}$ such that $\bint1{c_1c_2\dots c_n \dots c_{l_n-1}
s_n }\bigcap \Cset{V_k(\xi)} = \varnothing$. Hence, $S_{\xi}$ is a
nowhere dense set.

If  condition 1 (condition 2 resp.  condition 3) of the Theorem
holds, then, from the equality $S_\xi = \Cset{V_k(\xi)}$ and
Theorem \ref{thm:CVNk} (Theorem \ref{theorem4} resp. Theorem
\ref{theorem: <d}) it follows that $\lambda (S_{\xi})=0$.
\end{proof}

\textbf{Corollary.} If
\[
\prod_{k=1}^\infty \max_m\left\{p_{mk}\right\}=0
\]
and one of the conditions 1, 2, 3 of Theorem \ref{Contorovity}
holds, then the random variable $\xi$ has a Cantor-type singular
continuous distribution.

\textbf{Acknowledgement}

This work was supported by DFG 436 UKR 113/78, DFG 436 UKR 113/80,
SFB-611 projects and by Alexander von Humboldt Foundation. The last
three named authors gratefully acknowledge the hospitality of the
Institute for Applied Mathematics of the University of Bonn.

\end{document}